\documentclass [12pt]{amsart}

\usepackage{amsmath, accents} %

\usepackage{amssymb,amsxtra,amsfonts}
\usepackage{amsmath, accents} %
\usepackage{graphicx}

\usepackage{booktabs}   %
\usepackage[colorlinks]{hyperref}

\usepackage{xspace} 

\long\def\pb #1*/{}

\def\reE@DeclareMathSymbol#1#2#3#4{%
    \let#1=\undefined
    \DeclareMathSymbol{#1}{#2}{#3}{#4}}
\DeclareSymbolFont{symbolsC}{U}{txsyc}{m}{n}
\SetSymbolFont{symbolsC}{bold}{U}{txsyc}{bx}{n}
\DeclareFontSubstitution{U}{txsyc}{m}{n}
\reE@DeclareMathSymbol{\strictiff}{\mathrel}{symbolsC}{76}

\newcommand\beq{\begin{equation}}
\newcommand\eeq{\end{equation}}
\newcommand\bal{\begin{align*}}
\newcommand\eal{\end{align*}}   
\newcommand\bmx{\left(\begin{matrix}}
\newcommand\emx{\end{matrix}\right)}
\newcommand\bsmx{\left(\begin{smallmatrix}}
\newcommand\esmx{\end{smallmatrix}\right)}
\newcommand\bmxnp{\begin{matrix}}
\newcommand\emxnp{\end{matrix}}
\newcommand\bsmxnp{\begin{smallmatrix}}
\newcommand\esmxnp{\end{smallmatrix}}

\DeclareMathSymbol{\widehatsym}{\mathord}{largesymbols}{"62}
\newcommand\lowerwidehatsym{%
  \text{\smash{\raisebox{-1.3ex}{%
    $\widehatsym$}}}}
\newcommand\fixwidehat[1]{%
  \mathchoice
    {\accentset{\displaystyle\lowerwidehatsym}{#1}}
    {\accentset{\textstyle\lowerwidehatsym}{#1}}
    {\accentset{\scriptstyle\lowerwidehatsym}{#1}}
    {\accentset{\scriptscriptstyle\lowerwidehatsym}{#1}}
}

\newcommand{\wh}{\fixwidehat}

\newcommand{\usm}[1]{\underline{\smash{#1}}}

\newcommand{\bSi}{{\bf \Si}}

\newcommand{\onto}{\twoheadrightarrow}
\newcommand{\otno}{\twoheadleftarrow}  

\newcommand{\into}{\hookrightarrow}
\newcommand{\spq}{/\!\!/}

\newcommand{\st}{\ \bigl\vert\ }
\providecommand{\abs}[1]{\lvert#1\rvert}

\providecommand{\Euler}{\text{\rm Euler}} 
\providecommand{\MS}{\text{\rm MS}} 

\providecommand{\supp}{\text{\rm supp}}   

\providecommand{\floor}[1]{\lfloor#1\rfloor}
\providecommand{\ceiling}[1]{\lceil#1\rceil}

\def\part#1{\frac{\partial\phantom{q}}{\partial#1}}

\newcommand {\flb}{\lbrack\!\lbrack}
\newcommand {\frb}{\rbrack\!\rbrack}

\DeclareFontFamily{U}{wncy}{}
\DeclareFontShape{U}{wncy}{m}{n}{<->wncyr10}{}
\DeclareSymbolFont{mcy}{U}{wncy}{m}{n}
\DeclareMathSymbol{\Sh}{\mathord}{mcy}{"58} 

\usepackage{xparse}
\NewDocumentCommand{\dn}{e{_^}}{
  _{\IfValueT{#1}{#1}\vphantom{\smash[b]{|}}}
  ^{\IfValueT{#2}{#2}\vphantom{\smash[t]{|}}} }

\newcommand{\MB}{\mathcal{M}_{\text{\rm B}}}

\newcommand{\rank}{\mathop{\rm rank}}

\newcommand{\Prod}{\prod}

\DeclareMathOperator{\Hom}{Hom}         

\newcommand{\GL}{{\mathop{\rm GL}}}

\DeclareRobustCommand{\si}{\sigma}

\DeclareMathOperator{\Rep}{\rm Rep}

\newcommand{\diag}{{\mathop{\rm diag}}}

\newcommand{\hk}{{hyperk\"ahler }}   

\newcommand{\ba}{{\bf a}}

\newcommand{\bd}{{\bf d}}

\newcommand{\bl}{{\bf l}}
\newcommand{\bF}{{\bf F}}

\newcommand{\bQ}{{\bf Q}}

\DeclareSymbolFont{bbold}{U}{bbold}{m}{n}
\DeclareSymbolFontAlphabet{\mathbbold}{bbold}

\newcommand{\IA}{\mathbb{A}}
\newcommand{\IB}{\mathbb{B}}
\newcommand{\IC}{\mathbb{C}}

\newcommand{\IN}{\mathbb{N}}

\newcommand{\IP}{\mathbb{P}}                                     
                           
\newcommand{\IR}{\mathbb{R}}                           
\newcommand{\IS}{\mathbb{S}}

\newcommand{\IT}{\mathbb{T}}

\newcommand{\IV}{\mathbb{V}}

\newcommand{\IZ}{\mathbb{Z}}

\newcommand{\cB}{\mathcal{B}}

\newcommand{\cC}{\mathcal{C}}

\newcommand{\cM}{\mathcal{M}}

\newcommand{\cP}{\mathcal{P}}

\providecommand{\cR}{}
\renewcommand{\cR}{\mathcal{R}}

\newcommand{\cT}{\mathcal{T}}

\usepackage[mathscr]{euscript}
\newcommand{\sM}{       \mathscr{M}     }

\newcommand{\fM}{       \mathfrak{M}     }

\newcommand{\ga}{\gamma}

\newcommand{\Ga}{\Gamma}

\newcommand{\Si}{\Sigma}

\renewcommand{\th}{\theta}

\makeatletter
 \newlength{\typesize}
\makeatother

\newlength{\vvoff}
\newlength{\hhoff}

\newcommand{\pf}{\begin{bpf}}

\newcommand{\pfms}{\begin{bpfms}}
\newcommand{\epf}{\end{bpf}\hfill$\square$\\}           
\newcommand{\epfms}{\end{bpfms}\hfill$\square$\\}       

\newcommand{\idea}{\begin{bidea}}

\newcommand{\eidea}{\end{bidea}\hfill$\square$\\}           

\newcommand{\sk}{\begin{bsk}}    

\newcommand{\esk}{\end{bsk}\hfill$\square$\\}           
\newcommand{\sketch}{\begin{bsketch}}

\newcommand{\esketch}{\end{bsketch}\hfill$\square$\\}

\newtheorem {hypo}{\bf\hspace{-\parindent}Hypothesis}
\newtheorem {prop}[hypo]{Proposition}

\newtheorem {cor}[hypo]{Corollary}
\newtheorem {lem}[hypo]{Lemma}

\theoremstyle{definition}\newtheorem {defn}[hypo]{Definition}
\theoremstyle{definition}
\theoremstyle{definition} 
\theoremstyle{remark}\newtheorem{rmk}[hypo]{Remark}
\theoremstyle{remark}

\begin{document}

\title[Counting the fission trees]{Counting the fission trees and\\ nonabelian Hodge graphs (untwisted case)}
\author{Philip Boalch}%

\begin{abstract}
Any algebraic connection on a vector bundle 
on a smooth complex algebraic curve determines 
an irregular class and in turn a fission tree 
at each puncture. 
The fission trees 
are the discrete data classifying 
the admissible
deformation classes.
Here we explain how to count the 
fission trees with given slope and number of leaves, in the untwisted case. 
This also leads to a clearer picture of the ``periodic
table'' of the atoms
that play the role of building blocks 
in 2d gauge theory.

\end{abstract}

\maketitle

\section{Introduction}

The fission trees (\cite{twmcg} 3.18, \cite{rsode} p.29) 
encode the topology of 
the possible ways the structure group may be broken/fissioned 
by a meromorphic connection at a pole.
In the untwisted multiplicity $1$ case (our focus here), 
a fission tree  looks as follows:

\begin{figure}[ht]
    \centering
\includegraphics[width=0.8\textwidth]{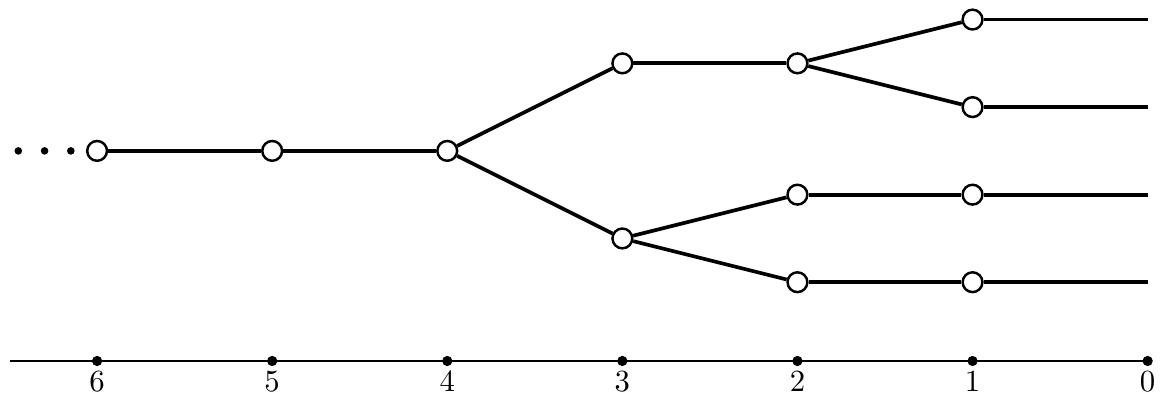}
  \caption{Example fission tree $\IT$ (with $4$ leaves and slope $3$).} 
    \label{fig: fissiontree}
\end{figure}

Thus it is connected, each node has an 
integer height, 
the leaves have height $0$, %
a node of 
height $1$ cannot be a branch node, and 
there is a highest branch node, the root, above 
which there is 
just one node for each integer.
The simplest numerical invariants of a fission tree
are the number $n$ of leaves 
and the slope $k$ ($=$ the height of the root minus $1$).

The aim of this note is to explain how to 
compute the number 
$\phi(k,n)$ of such fission trees of slope 
$k$ with $n$ leaves, up to isomorphism.
We will then adapt this to solve several closely related counting problems.

The above drawing of a fission tree is a bit redundant and it is often more convenient to draw a pruned form of the tree, obtained by truncating the tree at the root, and trimming the leaves off
(cutting at height $1$). Thus the pruned rendering of the above tree is as follows:

\begin{figure}[ht]
    \centering
\includegraphics[width=0.4\textwidth]{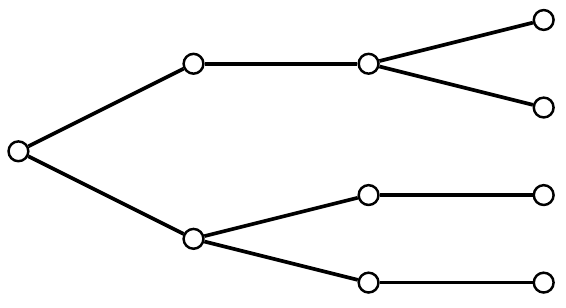}
    \label{fig: fissiontree2}
\end{figure}

For example it is not so hard to find the other eight fission trees with slope $3$ and $4$ leaves ($=$ the number of nodes of height $1$), as well as those with $2$ or $3$ leaves:

\begin{figure}[ht]
    \centering
\includegraphics[width=\textwidth]{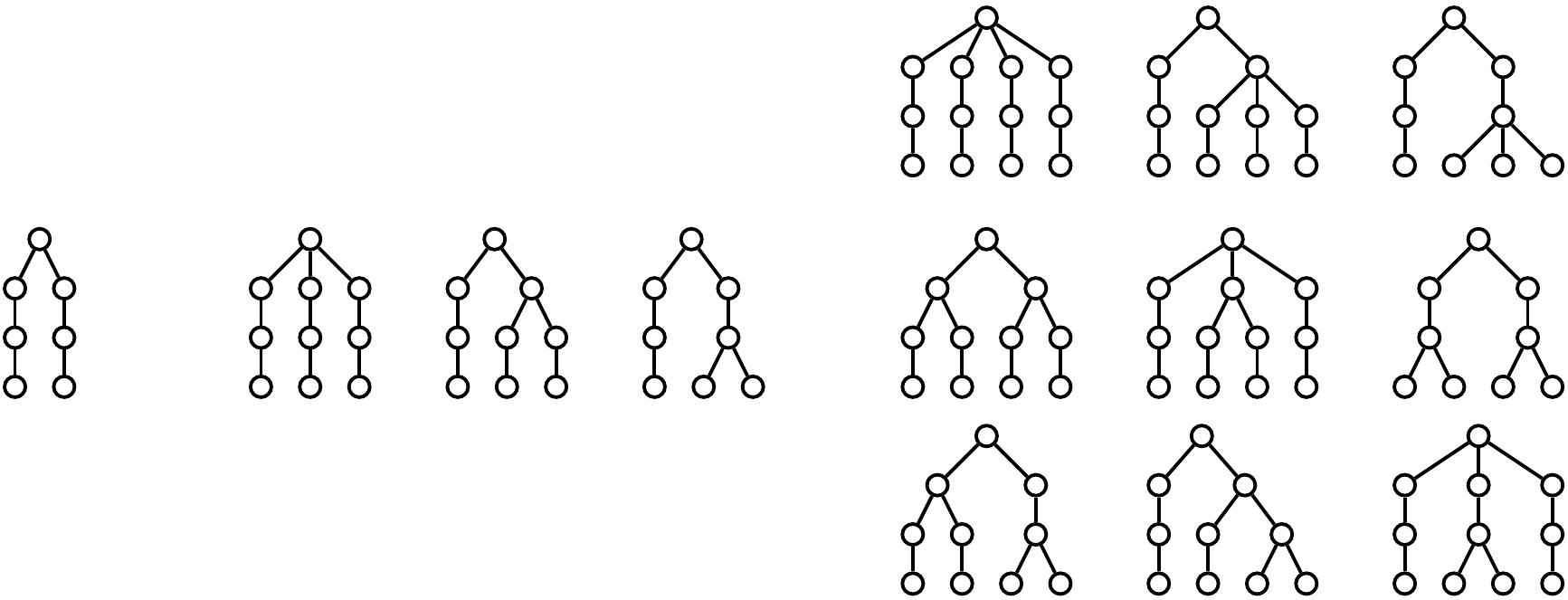}
    \caption{$1+3+9$  fission trees with slope $3$ and $2,3$ or $4$ leaves.}
    \label{fig: 139fissiontreesslope3}
\end{figure}

For other slopes/leaf counts, direct counting quickly becomes laborious.
We will explain how a method dating back to Euler
may be adapted to give 
a quick way to count the fission trees, for example 
leading to the construction of the following table:

\begin{table}[h]
 \begin{tabular}{| c || c | c | c | c | c | c | c | c | c | c | c | c | c | c   |} 
 \hline
 \phantom{$\Bigl\vert$}$k\setminus n$ &   1 &  2 &     3 &    4 &     5 &      6 &       7 &  8 &  9 &  10 \\  \hline\midrule
0 &   1 &  0 &     0 &    0 &     0 &      0 &       0 &  0 &    0 &  0 \\  
1 &   0 &  1 &     1 &    1 &     1 &      1 &       1 &  1 &    1 &  1 \\  
2 & 0 &  1 &  2 &     4 &    6 &    10 &     14 &      21 &  29 &  41  \\ 
3 & 0 &  1 &  3 &     9 &   20 &    47 &     96 &     201 &  394 &  775  \\ 
4 & 0 &  1 &  4 &    16 &   48 &   148 &    407 &    1121 &  2933 &  7612  \\ 
5 & 0 &  1 &  5 &    25 &   95 &   365 &   1271 &    4383 &  14479 &  47198  \\ 
6 & 0 &  1 &  6 &    36 &  166 &   766 &   3237 &   13466 &  53933 &  212645 \\ 
7 & 0 &  1 &  7 &    49 &  266 &  1435 &   7140 &   34853 &  164324 &  761829 \\ 
8 & 0 &  1 &  8 &    64 &  400 &  2472 &  14162 &   79430 &  431242 &  2301016 \\ 
9 & 0 &  1 &  9 &    81 &  573 &  3993 &  25893 &  164157 &  1009029 &  6094011 \\ 
10 & 0 &  1 &  10 &  100 &  790 &  6130 &  44392 &  314011 &  2156113 &  14544961\\ \bottomrule
\end{tabular}
\caption{Counting (untwisted, mult. 1) slope $k$ fission trees 
with $n$ leaves}
\end{table}

Note for example that row three: $0,1,3,9,20,47,96,...$ 
is not currently in the online encyclopaedia of integer sequences \href{https://oeis.org}{oeis.org}.

Three aspects of the underlying (somewhat elaborate) motivation are as follows. The reader  interested mainly in counting fission trees  might skip straight to \S2.

\subsection{Topological classification of irregular types} The main motivation for studying the fission trees comes from global Lie theory, or more precisely the classification 
of wild Riemann surfaces $\bSi$, and in turn the deformation classes of their character varieties $\MB(\bSi)$ and
nonabelian Hodge moduli spaces $\fM(\bSi,\cC)$
(i.e. we wish to understand the deformation classes of these complete 
hyperk\"ahler manifolds).
In simple terms one chooses a compact Riemann surface $\Si$, an integer $n\ge 1$ (the rank),  
a finite subset $\ba\subset \Si$ and an {\em irregular type} $Q_a$
at each point $a\in \ba\subset \Si$, and defines 
$\bSi$ to be the triple $(\Si,\ba,\bQ)$ where 
$\bQ=\{Q_a\st a\in \ba\}$.
This data 
determines an algebraic Poisson moduli space  $\MB(\bSi)$ 
as in \cite{gbs} Cor. 8.3 (and \cite{saqh,smid} in the generic case)
and its symplectic leaves have \hk metrics (they are the wild 
nonabelian Hodge moduli spaces $\fM(\bSi,\cC)$ constructed in \cite{wnabh}).
If we choose a local coordinate $z$ 
vanishing at $a\in \Si$ then an irregular type at $a$ is just an element of the form
$$Q_a=\frac{A_k}{z^k} + \cdots + \frac{A_2}{z^2}+ \frac{A_1}{z}$$
where the $A_i$ are all diagonal complex $n\times n$ 
matrices.
The choice of irregular types picks out a finite dimensional space of  meromorphic connections/Higgs bundles with good properties 
(generalising the tame case when $Q=0$).\footnote{
Given an irregular type $Q$, one considers meromorphic connections that are locally of the form $d-(dQ+ M(z) dz/z)$ for some holomorphic matrix $M$, as in the Hukuhara--Turrittin normal form, so that solutions of such connections involve the essentially singular functions $\exp(Q)$.
One can then look at moduli spaces of such connections, 
and their Betti descriptions \cite{JMU81,smid, saqh,gbs}
leading to the hyperkahler manifolds of \cite{wnabh}, 
classifying wild harmonic bundles on $\Si\setminus \ba$.}

Now, each irregular type $Q=\sum_1^k A_i/z^i$ determines a fission tree 
$\IT(Q)$, whose set of nodes of height $h\in \IN$ 
is the set of eigenspaces of the truncation
$$\tau_h(Q) = \sum_{i\ge h} A_i/z^i$$ of $Q$ (we take the truncation to be zero if $h> k$, so there is then just one node), as in \cite{twmcg} Lem. 3.22 (cf. also  \cite{rsode} Apx. C, \cite{gbs}).
For example $Q=A/z^3+B/z^2+C/z$ gives the tree in Fig. 1 
if $A=\diag(0,0,1,1), B=\diag(0,0,0,1), C=\diag(1,0,0,0)$.

The fission tree $\IT(Q)$ encodes the topology of the irregular type $Q$, and leads to the topological skeleton 
$\text{Sk}(\bSi):=(g,\bF)$ of $\bSi$, where $g$ is the genus of 
$\Si$ and $\bF=\sum_{a\in \ba}[\IT(Q_a)]$ is the fission forest, made up of the multiset of isomorphism classes of fission trees (as in \cite{twmcg} \S3.7). 
Now the key point here is that the moduli spaces $\MB(\bSi)$
vary smoothly (\cite{gbs} Thm. 10.2), provided that the 
wild Riemann surface  $\bSi=(\Si,\ba,\bQ)$ varies in an admissible fashion, i.e. that it varies smoothly and the topological skeleton $\text{Sk}(\bSi)$ is constant.
Thus the topological skeleta {\em parameterise} the admissible deformation classes.
In more detail
the  notion of admissible deformation of $\bSi$ was defined in \cite{gbs, twmcg} and in the current setting it says that:

$\bullet$ The underlying Riemann surface $\Si$ varies smoothly, and

$\bullet$ the points $\ba\subset \Si$ vary smoothly and do not coalesce, and

$\bullet$ each irregular type $Q_a$ varies smoothly and its fission tree $\IT(Q_a)$ does not change.

Thus it is important to classify the fission trees in order to understand the deformation classes of the wild character varieties $\MB(\bSi)$ and in turn  wild nonabelian Hodge moduli spaces 
$\fM(\bSi,\cC)$ (which are isomorphic to symplectic leaves of 
$\MB(\bSi)$)---in general distinct fission trees will lead to distinct deformation classes of moduli spaces.

Note that the underlying deformation spaces themselves, the moduli spaces $\sM_{g,\bF}$ of wild Riemann surfaces with given skeleton, also form interesting geometric objects, and their fundamental groups are the wild mapping class groups, \cite{p12} \S8, 
\cite{gbs, twmcg, DRTmwRS}. 
By definition a map $\IB\to \sM_{g,\bF}$ is the same thing as an admissible family over $\IB$ of wild Riemann surfaces (with skeleton 
$(g,\bF)$), as defined in \cite{gbs} \S10.
The spaces  $\sM_{g,\bF}$ generalise the usual moduli spaces of curves with marked points (that appear in the tame case, in effect when the fission trees all just have one leaf).
Enumerating the fission trees is the key step in enumerating the forests $\bF$ and thus the spaces $\sM_{g,\bF}$.

\subsection{Nonabelian Hodge graphs} \label{ssn: nabh graphs}
Secondly the fission trees are key to seeing which (non-affine) 
Kac--Moody 
root systems appear in 2d gauge theory. 
For example the nonabelian Hodge spaces of complex dimension two (the H3 surfaces) are related to ADE surface singularities 
(cf. \cite{saitoumemura, quad}, \cite{hit70} \S4.1), and  are classified by certain, special, affine Dynkin diagrams.
For example in the six Painlev\'e cases one gets the diagrams: 

\begin{figure}[ht]
    \centering
\includegraphics[width=0.9\textwidth]{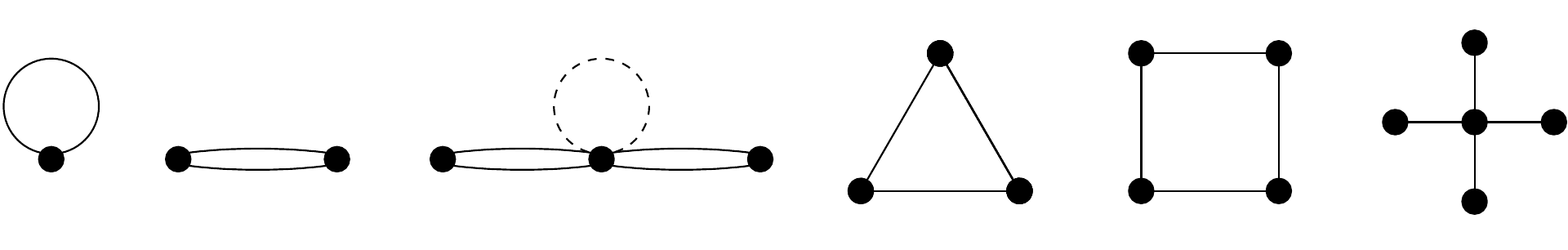}
\caption{The diagrams of the Painlev\'e moduli spaces 
1,2,3,4,5,6 (see \cite{diags}).}\label{fig: p123456}
 \end{figure}

In complex dimension four there are examples from hyperbolic Kac--Moody Dynkin graphs 
\cite{rsode} p.11-12, including the higher Painlev\'e systems 
\cite{slims} \S11.4, 
as well as several other examples \cite{kns4d, nakamura19}.  
This story generalises to 
 nonabelian Hodge spaces 
 of arbitrary dimension on the Riemann sphere:
 any such space 
has a diagram attached to it
 \cite{CB-additiveDS, CB-Shaw, rsode, slims, hi-ya-nslcase, yamakawa20, diags, doucot2021diagrams}
 (see \cite{hit70} \S1.6 for    
 more on nonabelian Hodge spaces). 
Some of the diagrams are just graphs 
(possibly with multiple edges, but no edge loops, 
or negative/dashed edges)
and others
 involve loops or negative edges (such as the dashed loop for Painlev\'e 3 above), that may be embedded in the untwisted 
case by pullback. 
In all cases one gets a generalisation of a Cartan matrix whose inner product controls the dimension of the moduli space \cite{diags, doucot2021diagrams}. 
It is natural to ask the basic question:
$$\boxed{\begin{matrix}\text{ What is the special class of graphs that appear as} \\
 \text{diagrams of nonabelian Hodge spaces?}\end{matrix}}$$
This question looks to be open, but we know a large class of examples of such nonabelian Hodge graphs, the {\em supernova graphs}. 
The supernova graphs are simple extensions of the  fission graphs, 
and their classification reduces to the classification of fission trees.
The two basic statements are as follows:

$\bullet$\ There is a special class of graphs, 
the {\em fission graphs}, 
that are in bijection with the 
(untwisted, multiplicity $1$) 
fission trees of slope $>1$. 
If $\IT$ is a fission tree 
with leaves $I$ then the fission graph $\Ga(\IT)$ of $\IT$ has vertices $I$ 
and the number of edges %
between two vertices $i,j$ in $\Ga$  
is $h_{ij}-2$ where $h_{ij}$ is the height of the
nearest common ancestor of the leaves $i,j$.
(The  slope $>1$ condition implies the graph is connected.)
They are very special: for example there are only 
$\phi(2,10)=41$
simple fission graphs with $10$ nodes, out of the  
11.7 million or so 
simple  graphs with $10$ nodes 
\href{https://oeis.org/A001349}{oeis/A001349}. 
For example the fission graphs of the nine fission trees of slope $3$ with $4$ leaves are:
\begin{figure}[ht]
    \centering
\includegraphics[width=.5\textwidth]{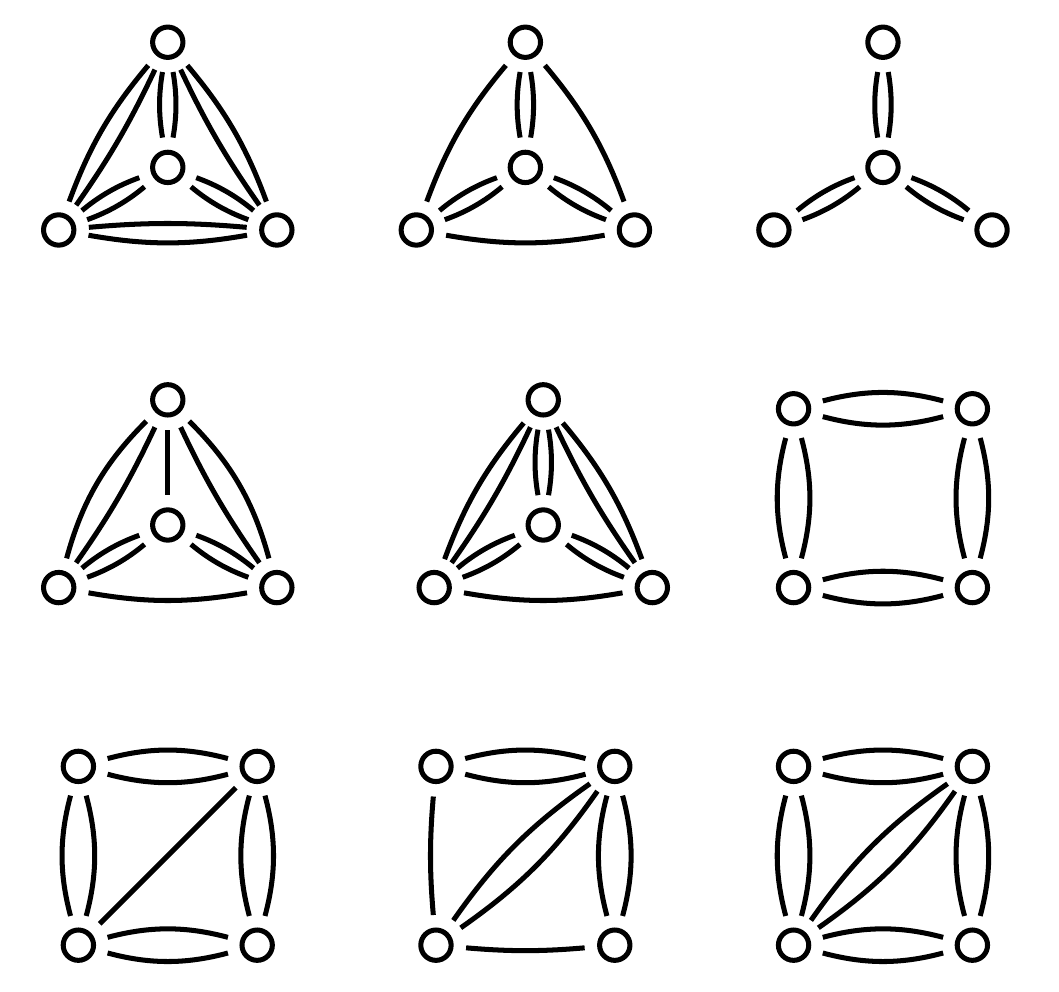}
    \caption{Fission graphs with $4$ nodes
    and maximal edge multiplicity $2$.}
    \label{fig: 9 fission graphs 34}
\end{figure}

$\bullet$\ Suppose $\Ga$ is a fission graph with nodes $I$
and we choose an integer $\ge 0$ 
for each node, i.e. we choose  
${\bf l}=(l_i)_{i\in I}\in \IN^I$.
Then the {\em supernova graph} $\wh \Ga(\bl)$  \label{defn: supernova}
determined by $(\Ga,\bl)$ 
is the graph obtained by gluing a leg of length $l_i$ 
onto the node $i\in I$ of $\Ga$ for each $i\in I$,
(\cite{slims} Def. 9.1, \cite{rsode} Apx. C).
The {\em core} of $\wh \Ga(\bl)$ is $\Ga$.
See the examples in Figs. \ref{fig: supernova}, \ref{345supenovagraphs}.

It is known  that
any supernova graph  $\wh \Ga(\bl)$ is a nonabelian Hodge graph\footnote{After seeing this paper 
J. Dou\c{c}ot showed the converse is not true \cite{doucotemail.2.25}.}.
In other words, there is a (nonempty) wild nonabelian Hodge moduli space 
$\fM(\bSi,\cC)$ with diagram   $\wh \Ga(\bl)$
(for some choice of marking of the formal monodromy orbits).
This follows 
(as explained in the simply-laced case in \cite{rsode} p.21)
using the Kac--Moody root system of the graph 
$\wh \Ga(\bl)$ and results of Vinberg \cite{vinberg1971} and Crawley--Boevey \cite{CB-mmap}---in the general non-simply-laced case, this uses Hiroe--Yamakawa's proof \cite{hi-ya-nslcase,  yamakawa20} of the quiver modularity conjecture of \cite{rsode} Apx. C
(giving a dictionary identifying Nakajima quiver varieties of supernova graphs as moduli spaces of meromorphic connections).

This gives further motivation for enumerating the fission trees (and thus the fission graphs). 
In particular this shows where these special types 
of Kac--Moody root systems appear in 2d gauge theory.

\begin{figure}[ht]
    \centering
\includegraphics[width=0.5\textwidth]{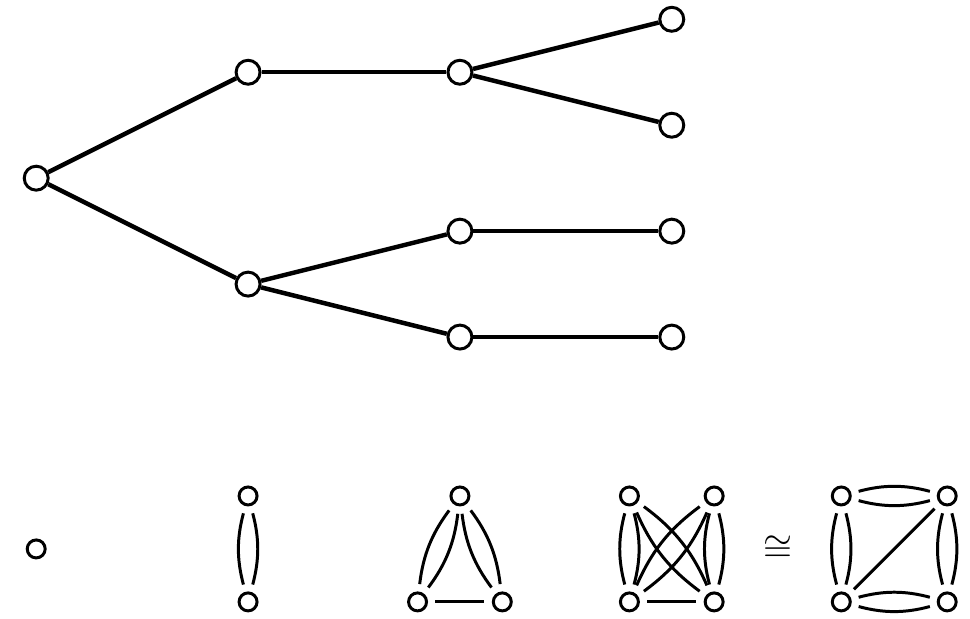}
  \caption{Fission construction of the fission graph of the fission tree of Fig. 1 (as drawn in \cite{rims15slides} p.79).}
    \label{fig: fissiontreeand graph}
\end{figure}

The original definition of the fission graph of an irregular type (\cite{rsode} Apx. C) involved a sequence of fission operations, as illustrated in Fig. \ref{fig: fissiontreeand graph}.

\begin{figure}[ht]
    \centering
\includegraphics[width=0.4\textwidth]{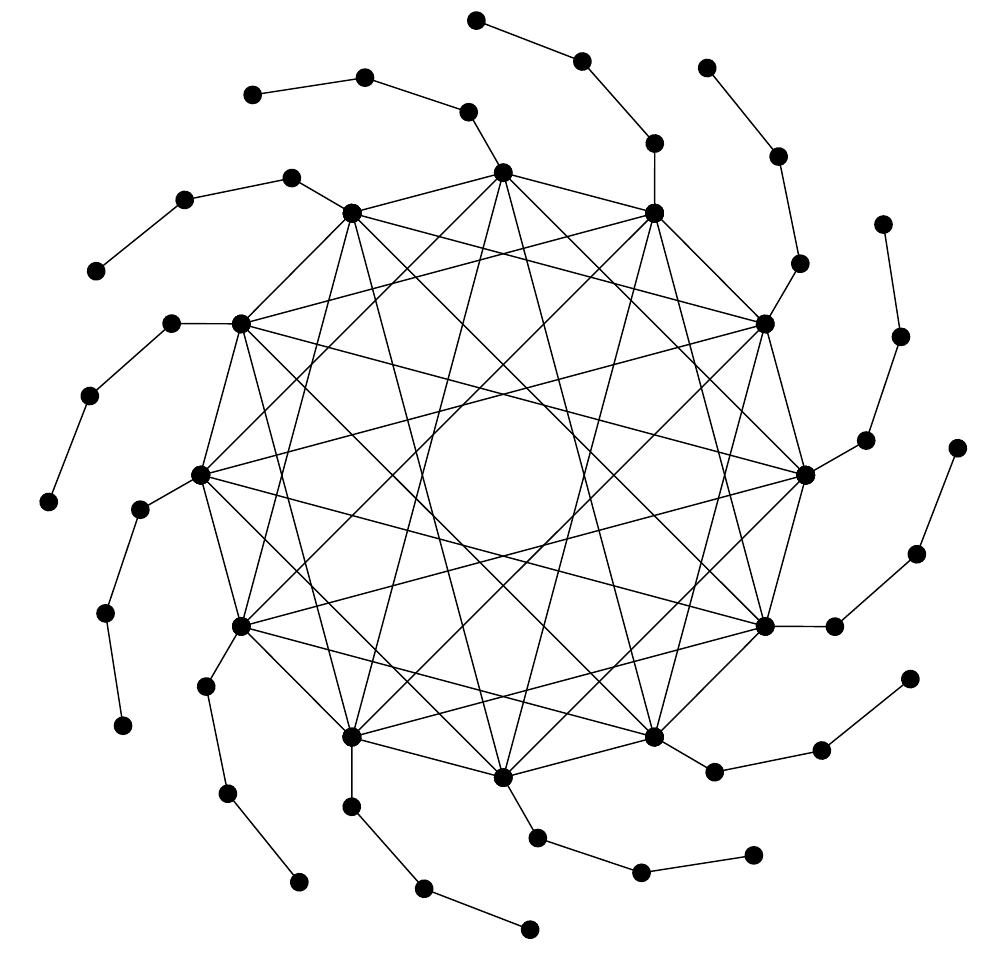}
\caption{Example supernova graph, with $l_i=3$ for all core nodes 
$i$.}\label{fig: supernova}
 \end{figure}

\begin{rmk}
The direct meaning  of the fission graph/tree in terms of Stokes arrows was pointed out in \cite{diags} p.65 (recall that the Stokes arrows \cite{tops} \S5.5 parameterise all the possible nontrivial  entries of the Stokes automorphisms, 
rephrasing the original definition in \cite{MR91}):
If we increase each edge multiplicity of a fission graph 
$\Ga(\IT)$ by one, then we get the {\em Stokes quiver} of the  fission tree $\IT$ (thus the Stokes quivers are also in bijection with the fission trees)---If $\IT=\IT(Q)$ is the fission tree of an irregular type, then the Stokes quiver contains all the Stokes arrows at the pole (the union of all the Stokes arrows for all the singular directions of $Q$, as in \cite{tops} \S5.5). 
See the examples of this passage from 
Stokes arrows/quivers to fission graphs 
on pp.80-83 of these 
\href{https://web.archive.org/web/20220221111633/https://webusers.imj-prg.fr/~philip.boalch/files/boalch_2022_diagrams,%20fission%20spaces%20and%20global%20Lie%20theory.pdf}{slides}.
\end{rmk}

\subsection{New multiplicative quiver varieties}

The ``quasi-Hamiltonian'' theory of Lie group valued moment maps \cite{AMM} 
gives a precise mathematical language for constructing symplectic moduli spaces 
via TQFT-type gluing of real surfaces.
Although originally set-up for compact Lie groups it may be transposed to the 
complex algebraic setting, and then leads to many new algebraic symplectic and Poisson moduli spaces \cite{saqh,vdb-doublepoisson,yamakawa-mpa, fission, gbs}.

The (new) multiplicative quiver varieties appear as follows:
Given a fission graph $\Ga$ with nodes $I$ and a dimension 
vector $\bd\in \IN^I$, one may define a Zariski open subset:
$$\cB(\Ga,\bd)=\Rep^*(\Ga,\bd) \subset \Rep(\Ga,\bd)$$
of the space $\Rep(\Ga,\bd)$ 
of representations of the graph $\Ga$ on the vector space $\IC^\bd$, 
the reduced fission space/invertible graph representations
(here we identify $\Ga$ as a quiver, by viewing each edge as two opposite arrows).
The reduced fission space 
is a quasi-Hamiltonian variety equipped with a moment map
$$\mu : \cB(\Ga,\bd) \to \GL_\bd(\IC):= \Prod_{i\in I}\GL_{d_i}(\IC).$$
In brief if $\Ga$ is the fission graph of 
a fission tree $\IT(Q)$ of an irregular type $Q$ at 
$\infty\in \IP^1$ then $\cB$ is 
isomorphic to the wild representation variety 
$\cR=\Hom_\IS(\Pi,G)$ of $\bSi=(\IP^1,\infty, Q)$, 
from \cite{gbs} Thm 1.1 
(and independent of $Q$ by \cite{gbs} Thm 10.4).

Now suppose we have an arbitrary finite graph $\Ga$ (possibly with multiple edges)
equipped with a colouring map 
$\ga:\Ga\to C$
to some set $C$ of colours, such that each monochromatic subgraph
$\Ga(c):=\ga^{-1}(c)$ is a fission graph.
Then for each dimension 
vector $\bd\in \IN^I$ and parameters $q\in (\IC^*)^I$ (where $I$ is the set of nodes of $\Ga$), there is a symplectic variety called the (new) 
multiplicative quiver variety:
$$\cM(\Ga,q,\bd) = \Rep^*(\Ga,\bd)\spq_{\!\!q} \,\GL_\bd(\IC)
=\mu^{-1}(q)/\GL_\bd(\IC)$$
where $\Rep^*(\Ga,\bd)$ is defined by fusing together the 
spaces  
$\cB(\Ga(c), \bd(c))$ for each colour $c$.
The classical multiplicative quiver varieties \cite{CB-Shaw, vdb-doublepoisson, yamakawa-mpa} are the special case when each fission graph is a single edge $\circ-\circ$, corresponding to the (unique) fission tree of slope $2$ and two leaves.
The simply-laced case  (when each fission graph has no multiple edges, i.e.  
the corresponding fission trees have slope $2$)
is described in detail in \cite{cmqv}, and the general case follows the same pattern (noting Rmk 5.4 of \cite{cmqv}).
In general the classical multiplicative quiver variety 
$\cM_{cl}(\Ga,q,\bd)$ looks to be a proper open subset of $\cM(\Ga,q,\bd)$ (i.e. when we modify the colouring map so that each edge has a distinct colour, one obtains a slightly smaller symplectic variety).

In any case we see the fission graphs $\Ga(\IT)$  parameterise  the basic building blocks  $\cB(\Ga,\bd)$ here, giving more motivation for the classification of the underlying fission trees  $\IT$.
They form a kind of ``periodic table'' for the building blocks, or ``atoms''.
(Via \cite{wnabh} %
they parameterise harmonic bundles on $\IA^1$, and it is hoped  
$\cM(\Ga,q,\bd)$ still parameterises harmonic bundles; they are special cases of fission varieties, obtained by fusing/gluing such pieces (cf. \cite{fission}, and the examples pp.115-144 of these \href{https://web.archive.org/web/20211224113708/https://webusers.imj-prg.fr/~philip.boalch/files/boalch_2017_oxford.and.ihes2017.pdf}{slides})).

\begin{rmk}
The general definition of the fission trees 
of a 
meromorphic connection on a vector bundle on a Riemann surface 
 in \cite{twmcg} involves the choice of a dimension vector (i.e. the multiplicity map $n$ in \cite{twmcg} 3.4.1). 
Here we are just looking at the untwisted case 
(and so the definition is essentially the same as that in 
\cite{rsode} Apx. C), and we are forgetting the choice of 
dimension vector in the classification (as it is easy to add this afterwards: see \S\ref{sn: arb mult}).
The extension to the twisted case may be accomplished by considering automorphisms of fission trees, and will be discussed elsewhere. 
\end{rmk}

\section{Counting}\label{sn: counting}

\subsection{Encoding via surjective maps}
First observe that the combinatorics of an untwisted 
 fission tree $\IT$ (as in Fig. \ref{fig: fissiontree}) 
 can be encoded as follows.  
Observe that $\IT$ is determined by 
its sets $J_1,J_2,\ldots$ of nodes of each height $\ge 1$, together with
 the ``parent maps'', which are the surjective maps:
$$\pi_i:J_i\onto J_{i+1}$$
taking any node to its parent. 
Each of the sets $J_i$ is finite, and the maps converge in the sense that 
$\abs{J_i}=1$ for all 
sufficiently large $i$.
The number of leaves of $\IT$  equals $\abs{J_1}$ and 
the
{slope} of the tree is 
$\min\{i-1\st \abs{J_i}= 1\}$.  
Thus we can just as well consider sequences 
of surjective maps of finite sets:
\beq
\cdots J_{k+2} \cong J_{k+1}=\{*\}\otno J_k \otno \cdots \otno J_2 \otno J_1
\eeq
where all the maps to the left of $*$ are
isomorphisms. 
It follows that the number $\phi(k,n)$  is (also) the number of isomorphism classes of such sequences of finite sets and surjective maps, such that 
$\abs{J_1}=n$ and 
$\min\{\,i\st \abs{J_i}=1\}=k+1$.

Up to relabelling this is the way the fission trees appeared in \cite{rsode} Apx. C, where  $J_i$ was 
the set of 
simultaneous eigenspaces of the leading 
coefficients
$\{A_k,\ldots, A_i\}$ 
of an {irregular type} 
of the form:
$$Q=\frac{A_k}{z^k} + \cdots + \frac{A_2}{z^2}+ \frac{A_1}{z}$$
where the $A_i$ are all diagonal matrices.
Thus the fission trees encode how the 
eigenspaces of the truncation $\tau_h(Q)=\sum_{i\ge h} A_i/z^i$
are refined as $h$ decreases 
(geometrically this encodes the types of  
growth of the maps $\exp(Q)$ that 
occur in solutions of the corresponding connections).  

\subsection{Counting}
We will say that a fission tree is {\em generic} if it has just one branch node (so exactly one of the parent maps is not bijective).
There is just one fission tree with one leaf, and it has slope $0$, and moreover any tree with slope $0$ has just one leaf, so: 
$$\phi(0,1)=1 \quad\text{and}\quad 
\phi(k,1)=\phi(0,n)=0 \text{ if $k>0, n>1$}.$$

For slope $1$ there is just one fission tree with any 
number $n>1$ of leaves (a generic fission tree), so that: 
$$\phi(1,n)=1\quad\text{for any $n>1,$ and  $\phi(1,1)=0$ as above.}$$

For slope $2$, 
we have a surjective map $J_1\onto J_2$ 
of finite sets where $n=\abs{J_1}$ is the number of leaves.
Thus the map gives a partition of $n$, 
where the parts are labelled by the set $J_2$.
By definition $\abs{J_2}>1$ so we are counting the partitions of $n$ 
with at least two parts.
Since there is just one partition with just one part it follows that:
$$\phi(2,n)=p(n)-1$$
where $p(n)$ is the number of 
integer partitions of $n$.

For slope $3$ and $n=1,2,3,4,5$ leaves, we find directly that there are
 $0,1,3,9, 20$ fission trees, respectively. 
 The first few are as drawn in Fig. 
 \ref{fig: 139fissiontreesslope3}.  
 Now this sequence $0,1,3,9,20$ yields only $5$ matches at 
 \href{https://oeis.org}{oeis.org}, 
 none of which seem relevant. Thus it seems no-one has counted the fission trees before.

It turns out however that 
the desired numbers appear as the {\em difference} of 
some easily 
generated sequences, using a method going back to Euler.
In particular we will easily be able to tabulate the 
numbers $\phi(n,k)$ far beyond 
anything we could count by hand.

\subsection{Changing the question}

If we have a fission tree of slope $k>1$ 
and delete the root node (at height $k+1$) 
to disconnect it, 
then we get a multiset of fission trees of slope 
$<k$ (where we may have to 
prune the tops a little bit, as the new roots may be lower down).

This suggests we should change the question:
Let 
$\Phi(k,n)$ be the number of
isomorphism classes of fission trees 
of slope $\le k$ and $n$ leaves.
This is useful since $\Phi(k,n)=\sum_{r\le k}\phi(r,n)$ and so it is clear that:

\begin{lem}
$\phi(k,n) = \Phi(k,n)-\Phi(k-1,n).$
\end{lem}

In turn, to compute $\Phi(k,n)$
recall that the {\em Euler transform} of an integer  sequence
$a_1,a_2,a_3,\ldots$ 
is the integer sequence $(b_n)$ given by:
$$1+\sum_1^\infty b_nu^n = 
\Prod_1^\infty(1-u^m)^{-a_m}.$$

It converts the counting of connected trees
with a given property into the counts of forests 
of trees with the same property 
(see \cite{sloaneplouffe} p.20).
For completeness a self-contained 
proof is in Appx. \ref{sn: appx Euler transform} below.

\begin{lem}
1)  $\Phi(2,n)=p(n)$ is the number of integer partitions of $n$,

2) for any fixed $k$ the sequence 
$\Phi(k,*)$ is the Euler transform of the sequence 
$\Phi(k-1,*)$.
\end{lem}
\pf
For 1), as mentioned above, 
note that the map $J_1 \onto J_2$
partitions the set $J_1$ into $\abs{J_2}$ parts.

For 2) first recall that a forest is a disjoint union of (not necessarily distinct) trees, i.e. it is an (unordered) multiset of isomorphism classes of trees. 
Then note that 2) 
means exactly that $\Phi(k,n)$
is the number of fission forests, each with $n$ leaves in total, such that each 
tree in the forest has slope $\le k-1$.
Any such forest determines a unique 
fission tree with $n$ leaves and of slope $\le k$:
View each tree in the forest as a full (nontruncated) tree
and then truncate it at height $k+1$.
Then glue all the top nodes (height $k+1$) of all the trees together, 
to a single node $\{*\}$ (i.e. take the wedge sum).
This gives a fission tree of slope  $\le k$ (if the original forest had just one tree
then the resulting tree has slope $<k$, and needs to be pruned to get a truncated tree).
Conversely any fission tree counted by 
$\Phi(k,n)$ determines such a forest by deleting the root node (and pruning the tops).  
\epf

Note that 2) holds in general and the
numbers $\Phi(2,n)$ are the Euler transform 
of the sequence 
$[1,1,1,1,1,\ldots]$ essentially counting the 
fission trees of slope $1$.
This is the original instance:
$$1+\sum_1^\infty p(n)u^n = 
\Prod_1^\infty(1-u^m)^{-1}$$
of the Euler transform,
giving the count $p(n)=\Phi(2,n)$ of partitions.

Note (of course) that 
the case before this also holds:
$$1+u+u^2+u^3+\cdots = (1-u)^{-1} $$
which says that the all $1$'s sequence
$\Phi(1,n)$
is the Euler transform of the sequence 
$[1,0,0,0,\ldots]$ counting the fission trees of slope $0$.

The next case
$$1+\sum_1^\infty \Phi(3,n)u^n = 
\Prod_1^\infty(1-u^m)^{-p(m)}$$
is discussed by Cayley 1855 \cite{Cayley1855} p.316,
who was counting ``partitions of partitions'' 
(i.e. unordered multisets of partitions, whose total sum is $n$), 
which we can now identify with fission trees of slope $\le 3$.
They are sometimes called 
double partitions \cite{kaneiwa} or $2$-dimensional partitions\footnote{
But beware that a ``plane partition'' is a special type of 
$2d$ partition 
(and they have different counts for $n\ge 4$).}. 
Thus the number  $\phi(3,n)=\Phi(3,n)-p(n)$ we are interested in counts the {\em proper} double partitions of $n$ 
(i.e. those that do not just consist of a single partition).
As Cayley wrote, the list of counts $\Phi(3,n)$  is:
$$1,3,6,14,27,58,111,223,424,817,\ldots$$
which is \href{https://oeis.org/A001970}{oeis/A001970}, 
and $p(n)$ is:
$$1, 2, 3, 5, 7, 11, 15, 22, 30, 42, \ldots$$
so that the count $\phi(3,n)=\Phi(3,n)-p(n)$ 
of slope $3$ fission trees with $n$ leaves is:
$$0,1,3,9,20,47,96,201,394,775,\ldots$$
extending the list of counts we constructed above.  

For example there are $14$ double 
partitions of $4$, of which $5=p(4)$
are just single partitions. The
remaining $9$, the proper double partitions of $4$, are as follows:
$$\begin{matrix}
[[1][1][1][1]]& 
[[111][1]] &
[[3][1]]  \phantom{\Bigl\vert}\\
[[11][11]] &
[[11][1][1]] &
[[2][2]] \phantom{\Bigl\vert} \\
[[2][11]] &
[[21][1]]& 
[[2][1][1]] \phantom{\Bigl\vert}
\end{matrix}
$$
and correspond to the slope $3$ fission trees with $4$ leaves,  
drawn on the right of 
Fig.  \ref{fig: 139fissiontreesslope3}.
The corresponding fission graphs are 
in Fig. \ref{fig: 9 fission graphs 34}.

In any case it is now easy to compute the numbers $\phi(k,n)$ and construct Table 1.

\begin{table}[h]
 \begin{tabular}{| c || c | c | c | c | c | c | c | c | c | c | c | c | c | c   |} 
 \hline
 \phantom{$\Bigl\vert$}$k\setminus n$ &   1 &  2 &     3 &    4 &     5 &      6 &       7 &  8 &  9 &  10 \\  \hline\midrule
0 &   1 &  0 &     0 &    0 &     0 &      0 &       0 &  0 &    0 &  0 \\  
1 &   0 &  1 &     1 &    1 &     1 &      1 &       1 &  1 &    1 &  1 \\  
2 & 0 &  1 &  2 &     4 &    6 &    10 &     14 &      21 &  29 &  41  \\ 
3 & 0 &  1 &  3 &     9 &   20 &    47 &     96 &     201 &  394 &  775  \\ 
4 & 0 &  1 &  4 &    16 &   48 &   148 &    407 &    1121 &  2933 &  7612  \\ 
5 & 0 &  1 &  5 &    25 &   95 &   365 &   1271 &    4383 &  14479 &  47198  \\ 
6 & 0 &  1 &  6 &    36 &  166 &   766 &   3237 &   13466 &  53933 &  212645 \\ 
7 & 0 &  1 &  7 &    49 &  266 &  1435 &   7140 &   34853 &  164324 &  761829 \\ 
8 & 0 &  1 &  8 &    64 &  400 &  2472 &  14162 &   79430 &  431242 &  2301016 \\ 
9 & 0 &  1 &  9 &    81 &  573 &  3993 &  25893 &  164157 &  1009029 &  6094011 \\ 
10 & 0 &  1 &  10 &  100 &  790 &  6130 &  44392 &  314011 &  2156113 &  14544961\\ \bottomrule
\end{tabular}

Table 1. {Counting untwisted multiplicity 1 fission trees 
with slope $k$ and $n$ leaves}
\end{table}

One can compute much larger versions of this table with the following  maple code.

\tiny

\noindent\rule{12cm}{.05cm}

\vspace{-.2cm}\noindent\verb
# Maple code for n terms of kth row of table 1  

\vspace{-.2cm}\noindent\verb
# Euler transform adapted from https://oeis.org/transforms.txt

\vspace{-.2cm}\noindent\verb
k:=4: n:=12: 


\vspace{-.2cm}\noindent\verb
EULER:=proc(a)  local b,c,i,d: b:=[]:c:=[]:

\vspace{-.2cm}\noindent\verb
for i to nops(a) do c:=[op(c), add( d*max(0,1-irem(i,d))*a[d], d=1..i)]: od:

\vspace{-.2cm}\noindent\verb
for i to nops(a) do b:=[op(b),(1/i)*(c[i]+add( c[d]*b[i-d], d=1..i-1))]: od:

\vspace{-.2cm}\noindent\verb
RETURN(b);  end:

\vspace{-.2cm}\noindent\verb
(EULER@@(k-2))([1$n]):EULER(

\vspace{-.2cm}\noindent
\rule{12cm}{.05cm}

\normalsize

By inspection, there appear to be simple formulae for the
first few columns:
$$
\phi(k,3)=k,\quad 
\phi(k,4)=k^2,\quad 
\phi(k,5)=k(5k^2 - 3k + 4)/6,
$$
the last  of which equals \href{https://oeis.org/A203552}{oeis/A203552} 
(although no links to trees are noted there).
It seems one gets a (numerical) polynomial of degree 
$n-2$ for the $n$th column in general\footnote{We will prove that 
$\phi(k,4)=k^2$ in the appendix---the others should be viewed as experimental observations that one might want to 
try to prove directly.},
for example:
$\phi(k,6)=k(2k^3-2k^2+4k-1)/3$.

The simply-laced fission graphs (those with no multiple edges) are exactly the complete multipartite graphs, 
and they are parameterised by 
the {\em nontrivial} integer partitions 
(i.e. partitions with at least two parts),
coming from the slope $2$ fission trees. 
Recall that the complete $k$ partite 
graph with $n$ nodes determined by 
the partition $n_1+\cdots+n_k =n$ of $n$, is the 
graph $\Ga(n_1,\ldots,n_k)$ with $n$ nodes partitioned 
into parts of sizes $n_1,\ldots,n_k$ such that 
two nodes are connected by an edge if and only if they
are in different parts.
For example the core of the supernova graph in Fig. 
\ref{fig: supernova} is the complete bipartite graph $\Ga(6,6)$.
Looking at the slope $2$ row of the table,
we see the number of simply-laced fission graphs 
with $\le 6$ nodes is $1+2+4+6+10=23$. 
This includes the $5$ star-shaped graphs
$\Ga(1,n),n=1,2,3,4,5$. 
The remaining $18=1+3+5+9$ were drawn 
in \cite{rsode} Fig. 1 (and in \cite{slims} Fig. 3):

\begin{figure}[ht]
    \centering
\includegraphics[width=0.7\textwidth]{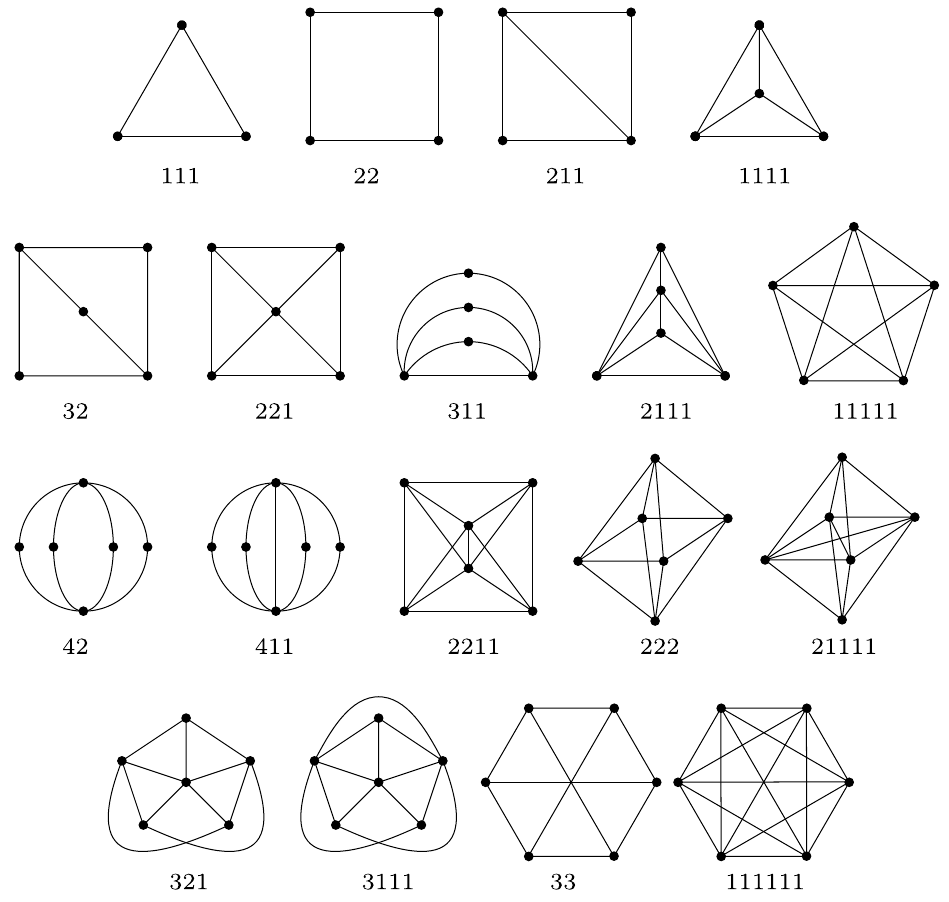}
\caption{The $18$ non-starshaped simply-laced 
  fission graphs with $\le 6$ nodes}\label{fig: simplefissiongraphs}
 \end{figure}

\section{Allowing arbitrary multiplicities}\label{sn: arb mult}

In the above discussion we counted the 
number of untwisted fission trees 
that appear when we forget the 
multiplicity of the leaves
(in effect setting the multiplicity of each leaf to be $1$).
Here we will explain how to add back in the multiplicities, 
and find 
essentially the same table but shifted by $1$.  
(Recall that the general definition of fission trees in \cite{twmcg} Defn. 3.18
includes the data of a multiplicity map 
$\underline{n}:\IV_0\to \IZ_{\ge 1}$ 
assigning an integer multiplicity to each leaf $\IV_0$ of the tree.)

Define the {\em rank} of a fission tree $\IT$ to be the sum of all the multiplicities:
$$\rank(\IT)=\sum_{i\in \IV_0} \underline{n}(i)$$
where $\IV_0$ is the set of leaves and $\underline{n}$ is the multiplicity map.
(If $\IT=\IT(Q)$ comes from an $n\times n$ irregular type $Q$, then $\IT$ has rank $n$
and the multplicities come from the multiplicities of the diagonal entries of $Q$.) 
Of course in the multiplicity $1$ setting above, the rank is the same as the number of leaves.
Also, since $\IV_1\cong \IV_0$ in the present untwisted setting, 
and we like to prune the leaves, we can just as well view 
$\underline{n}$ as a map 
$\underline{n}:\IV_1\to \IZ_{\ge 1}$.

Define $\psi(k,n)$ to be the size of the set of isomorphism classes of 
all untwisted fission trees of slope $k$ and rank $n$,
with arbitrary multiplicities. 
Recall that $\phi(k,n)$ was the size of the set of isomorphism classes of 
untwisted multiplicity $1$ fission trees of slope $k$ and rank $n$.
Clearly  $\psi(0,n)=1$ for all $n$ as there is just one fission tree with slope zero 
and rank $n$ 
(the tree with one leaf of multiplicity $n$).
For higher slopes $k$ one can reduce to the case already studied:

\begin{prop}
If $k\ge 1$ then there is a bijection between the set of isomorphism classes of 
untwisted fission trees of slope $k$ and rank $n$,
and the set of isomorphism classes of 
multiplicity $1$  untwisted fission trees of slope $k+1$ and rank $n$, so that  $\psi(k,n)=\phi(k+1,n)$ for any $n$.
\end{prop}
\pf
The bijection arises as follows: Suppose $\IT$ is  a multiplicity $1$  untwisted fission tree of slope $k+1\ge 2$ 
and rank $n$, with nodes $\IV$.
Define a tree $\IT'$  with nodes $\IV_i':=\IV_{i+1}$ for $i\ge 1$ (and the same edges between these nodes), 
in effect chopping off all of  $\IT$ below  height $2$, and shifting the height by $1$.
Then define the multiplicity 
$\underline{n}(i)$ of any element $i\in \IV_1'$
to be the number $\ge 1$ of children of the node $i\in \IV_2\subset \IT$. 
Thus $\IT'$ is an untwisted fission tree of slope $k$ and rank $n$,
(with possibly nontrivial multiplicities). 
This process is clearly bijective.
\epf

\begin{figure}[ht]
    \centering
\includegraphics[width=.45\textwidth]{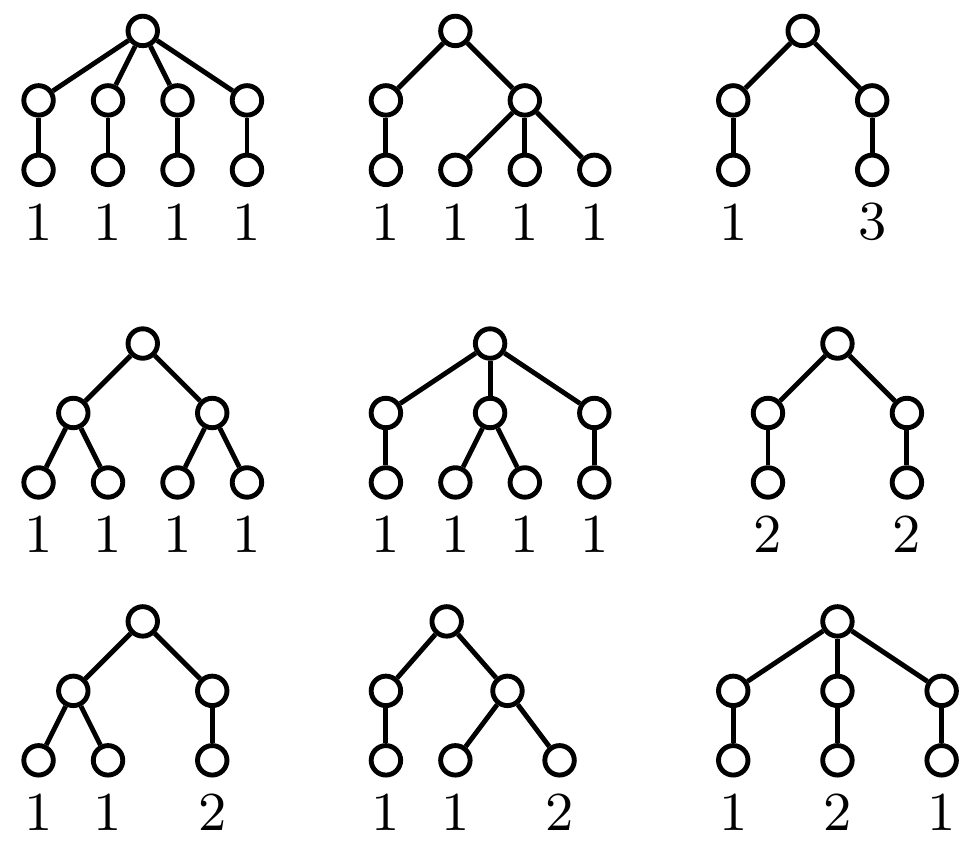}
    \caption{The $\psi(2,4)=9$ untwisted fission trees 
    of slope $2$, rank $4$.}
    \label{fig: 9fissiontreesslope2}
\end{figure}

Thus it is easy to modify table 1 to get the counts of all untwisted fission trees.
\begin{table}[h]
 \begin{tabular}{| c || c | c | c | c | c | c | c | c | c | c | c | c | c | c   |} 
 \hline
 \phantom{$\Bigl\vert$}$k\setminus n$ &   1 &  2 &     3 &    4 &     5 &      6 &       7 &  8 &  9 &  10 \\  \hline\midrule
0 &   1 &  1 &     1 &    1 &     1 &      1 &       1 &  1 &    1 &  1 \\  
1 & 0 &  1 &  2 &     4 &    6 &    10 &     14 &      21 &  29 &  41  \\ 
2 & 0 &  1 &  3 &     9 &   20 &    47 &     96 &     201 &  394 &  775  \\ 
3 & 0 &  1 &  4 &    16 &   48 &   148 &    407 &    1121 &  2933 &  7612  \\ 
4 & 0 &  1 &  5 &    25 &   95 &   365 &   1271 &    4383 &  14479 &  47198  \\ 
5 & 0 &  1 &  6 &    36 &  166 &   766 &   3237 &   13466 &  53933 &  212645 \\ 
6 & 0 &  1 &  7 &    49 &  266 &  1435 &   7140 &   34853 &  164324 &  761829 \\ 
7 & 0 &  1 &  8 &    64 &  400 &  2472 &  14162 &   79430 &  431242 &  2301016 \\ 
8 & 0 &  1 &  9 &    81 &  573 &  3993 &  25893 &  164157 &  1009029 &  6094011 \\ 
9 & 0 &  1 &  10 &  100 &  790 &  6130 &  44392 &  314011 &  2156113 &  14544961\\ 
10& 0& 1& 11& 121& 1056& 9031& 72248& 564201& 4280870& 31910879\\ \bottomrule
\end{tabular}

Table 2. {$\psi(k,n)$ counting all untwisted fission trees 
with slope $k$ and rank $n$}
\end{table}

\begin{rmk}
Whereas a fission tree
includes a multiplicity map (assigning an integer to each leaf),
by definition a quiver or graph
does not include a dimension vector. 
Thus since $\phi(k,n)$ counts the fission graphs 
with $n$ nodes and maximal edge multiplicity exactly $k-1$, 
the number $\psi(k,n)$ counts the ``equipped fission graphs'' 
$(\Ga,\bd)$
(i.e.  $\Ga$
equipped with a dimension vector $\bd=(d_i), d_i\ge 1$) such that the maximal edge multiplicity is exactly $k-1$ and the rank $(\sum d_i)$ is 
$n$.  
Observe it is actually the equipped fission graphs that parameterise the reduced fission spaces $\cB(\Ga,\bd)$.
\end{rmk}

\begin{figure}[h]
    \centering
\includegraphics[width=.45\textwidth]{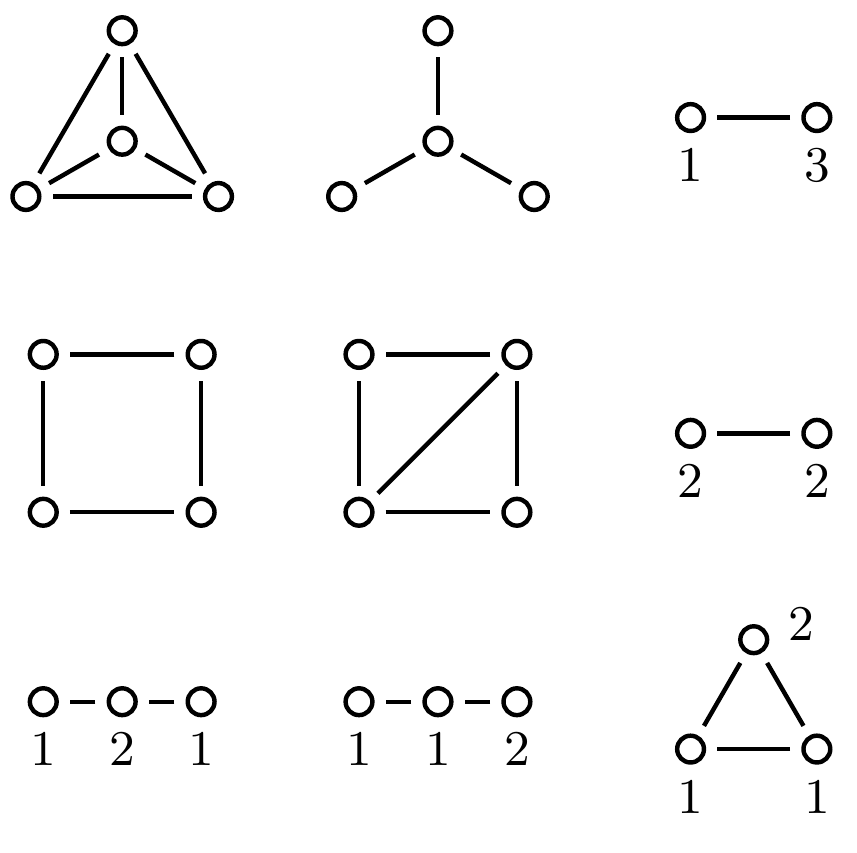}
    \caption{The $\psi(2,4)=9$ equipped fission graphs of rank $4$ with maximum edge multiplicity exactly $1$ (unmarked nodes have 
    dim. $1$).}
    \label{fig: 9equippedfissiontreesslope2}
\end{figure}

\section{Counting the supernovas}\label{sn: counting supernovas}

Recall from \S\ref{ssn: nabh graphs}  that a supernova graph is a graph of the form $\wh \Ga(\bl)$
obtained by gluing a leg of length $l_i\ge 0$ on to each node $i$ of a fission graph $\Ga$
(and that they seem to most closely  play 
the role of Dynkin graphs for moduli spaces
 in 2d gauge theory). 
Define $\si(k,n)$ to be the number of 
supernova graphs with 
$n$ nodes and maximal edge multiplicity $k$. 
In the last section we counted the ``equipped fission graphs'', 
i.e. the fission graphs $\Ga$ equipped with a dimension vector $\bd$ assigning an integer $d_i\ge 1$ to each node.  
Thus we have a map from equipped fission graphs to supernova graphs, defined by 
taking $l_i=d_i-1$. 
The number of nodes of the resulting 
supernova graph equals the rank $\sum d_i$ 
of the equipped fission graph.

This map is not always bijective:
For $k=1$ (the simply-laced supernova graphs, i.e. those that are simple graphs), 
the core is not uniquely determined, and so there is overcounting:
The issue is clear already in Fig. \ref{fig: 9equippedfissiontreesslope2}:
the $9$ equipped fission graphs there only determine $6$ distinct supernova graphs ($D_4$ appears twice and $A_4$ three times).

However this map gives a bijection most of the time, and thus enables us to count the supernova graphs, provided that the core is not star-shaped.

\begin{prop}
Suppose $k\ge 2$.
The map above gives a bijection between the 
equipped fission graphs
of rank $n$ with maximal edge multiplicity exactly $k$, 
and the supernova graphs with $n$ nodes and maximal edge multiplicity $k$.
Thus $\si(k,n)=\psi(k+1,n)=\phi(k+2,n)$.
\end{prop}
\pf
By definition all supernova graphs occur in this way, so the map is surjective.
The core of such a supernova graph is uniquely 
determined (exercise), and so distinct equipped fission graphs determine distinct supernova graphs.
\epf

For $k=1$  the issue is just with star-shaped core graphs 
(a supernova graph that is not star-shaped has a unique core) and we can count them directly to correct the counting:
$$\si(1,n)=\psi(2,n) - N_1(n) + N_2(n)$$
where
\begin{align*}
N_1(n) &:=\#\{\text{equipped, star-shaped, rank $n$ fission graphs}\}\\
N_2(n) &:=\#\{\text{star-shaped (supernova) graphs with $n$ nodes}\}.
\end{align*}

For example for $n=4$ as in 
Fig. \ref{fig: 9equippedfissiontreesslope2},
we have $N_1=5, N_2=2$ so that 
there are $\si(1,4)=9-5+2=6$ simple supernova graphs with $4$ nodes.

\begin{lem}
$N_1(n)= \floor{n/2}+ p(2) +\cdots +p(n-1) +2-n$ where $p(k)$ is the partition function.
\end{lem}
\pf
This is the same as counting the 
(untwisted mult. 1) fission trees of slope $2$ with $n$ leaves such that: 
a) exactly two nodes have height $3$ and 
b) at most one of the nodes at height $3$ has $>1$ child.  
If both nodes at height $3$ have $1$ child, then we are counting the partitions of $n$ with $2$ parts: there are $\floor{n/2}$
 of them.
If the first node at height $3$ has $1$ child, and the second has 
$>1$ child, then they have $n_1+n_2=n$ total descendants respectively, then the count is given by the number of  partitions of  $n_2$ into $\ge 2$ parts.  
Thus we get 
$$(p(n-1)-1)  + (p(n-2)-1) +\cdots +(p(2)-1)$$
going through the cases when $n_1=1,2,\ldots,n-2$.
Adding these up gives the answer.
\epf

If we define 
$\th(n)=\sum_0^n p(k)=2+\sum_2^np(k)$ 
(as in the sequence 
 \href{https://oeis.org/A000070}{oeis/A000070})
 then the lemma 
 says that $N_1(n)= \th(n-1) - \ceiling{n/2}$.

Next we can count the stars with $n$ nodes.
(Here a star-shaped graph is a finite connected simple graph with at most one vertex of degree $>2$.)

\begin{lem}
There are $N_2(n)=p(n-1)-\floor{(n-1)/2}$ star-shaped graphs with $n$ nodes, where $p(n)$ is the number of partitions of $n$.
\end{lem}

\pf
If there are less than $3$ legs then it is just a type $A_n$ Dynkin graph.
Otherwise we need to specify the length of the $\ge 3$ legs, i.e. choose a partition of $n-1$ with at least $3$ parts.
There is one partition with $1$ part and $\floor{(n-1)/2}$ with $2$ parts, and so we get 
$p(n-1)-1-\floor{(n-1)/2}$ partitions of $n-1$ with at least $3$ parts.
Adding on $1$ corresponding to $A_n$ gives the answer.
\epf

\begin{cor}
The number of simple supernova graphs with $n$ nodes is: 
$\si(1,n)=\psi(2,n) + 1 - \th(n-2) $
where $\th(k)=\sum_0^k p(i)$.
\end{cor}
\pf
This is immediate now, since 
$\th(n-1) - p(n-1)=\th(n-2)$ 
and $\ceiling{n/2}-\floor{(n-1)/2}=1$.
\epf

This is now easily computed 
(recall that $\psi(2,n)= \Phi(3,n)-p(n)$ where 
$\Phi(3,n)$ is the double partition function
\href{https://oeis.org/A001970}{oeis/A001970}). 
In summary the 
supernova counts are as in Table 3 below  
(which continues as in Table 2,
with $k$ shifted by $1$, as 
 $\si(k,n)=\psi(k+1,n)$ for $k\ge 2$).

\begin{figure}[ht]
    \centering
\includegraphics[width=\textwidth]{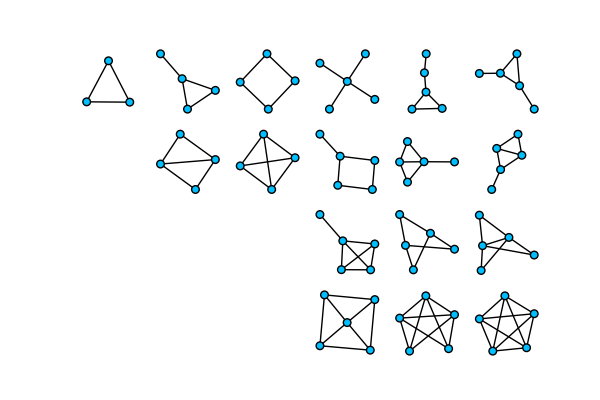}
     \caption{The non-Dynkin simple supernova graphs with  $3, 4$ or $5$ nodes, 
     including  the graphs  of the simply-laced H3 surfaces $\widehat A_2, \widehat A_3, \widehat D_4$. }
    \label{345supenovagraphs}
\end{figure}

\begin{table}[h]
 \begin{tabular}{| c || c | c | c | c | c | c | c | c | c | c | c | c | c | c   |} 
 \hline
 \phantom{$\Bigl\vert$}$k\setminus n$ &   1 &  2 &     3 &    4 &     5 &      6 &       7 &  8 &  9 &  10 \\  \hline\midrule
1 & 0 &  1 &  2 &     6 &   14 &    36 &     78 &     172 &  350&  709  \\ 
2 & 0 &  1 &  4 &    16 &   48 &   148 &    407 &    1121 &  2933 &  7612  \\ 
3 & 0 &  1 &  5 &    25 &   95 &   365 &   1271 &    4383 &  14479 &  47198  \\ 
\vdots & \vdots & \vdots & \vdots & \vdots & \vdots & \vdots & \vdots & \vdots & \vdots & \vdots\\
\bottomrule
\end{tabular}

Table 3. {$\si(k,n)$: supernova graphs 
with $n$ nodes \& max. edge mult. $k$.}
\end{table}
In fact we are most interested in the non-Dynkin supernova graphs, 
since the
Dynkin cases have no imaginary roots so only lead to rigid 
irreducible connections 
(nonabelian Hodge spaces of dimension zero).
Thus we can subtract the counts %
of the ADE Dynkin graphs 
$ A_1, A_2, A_3, (A_4,D_4) , (A_5,D_5), (A_6,D_6,E_6), \ldots$, 
to obtain Table 4 below.
The last column  continues as $\sigma(1,n)-2$ for $n>8$, counting the non-Dynkin simply-laced supernova graphs. 
See Fig. \ref{345supenovagraphs} for the first few examples.

\setcounter{table}{3}
\begin{table}[ht]
\begin{center}
\caption{Counting the simply-laced supernova graphs.}
 \begin{tabular}{|c | c | c || c|| c |} 
 \hline
  \#\text{nodes} & \#\text{starshaped}\phantom{$\Bigl\vert$} &  \#\text{other} &  Total & Total non-Dynkin\\ [0.5ex] 
 \hline\hline
2 &   1 &    0 &    1 & 0\\ \hline
3&  1&  1&  2&  1\\ \hline
4&  2&  4&  6&  4\\ \hline
5&  3&  11&  14&  12\\ \hline
6&  5&  31&  36&  33\\ \hline
7&  8&  70&  78&  75\\ \hline
8&  12&  160&  172&  169\\ \hline
9&  18&  332&  350&  348\\ \hline
10&  26&  683&  709&  707\\ \hline
 \hline
\end{tabular}
\end{center}
\end{table}

\section{Extended fission trees}

The fission trees really only capture {\em part} 
of the breaking of structure group at a pole, namely that due to the irregular/wild part  of the connection.
In general there is further breaking of the structure group 
of a slightly different nature ({\em tame fission}) from the tame/logarithmic part of the connection,
passing to the associated graded of the tame
(Levelt--Simpson) filtrations 
and then taking generalised eigenspaces
(cf. the tame fission
spaces $\mathbb{M}$ and the weighted conjugacy classes $\wh\cC$ in \cite{logahoric} \S4, and the review of the wild
nonabelian Hodge boundary data in \cite{ihptalk}).  
In terms of deformation classes, this amounts to upgrading the fission trees as follows (adding the ``fine hidden structure'' of the leaves).

\begin{defn}
A \usm{tame fission tree} $\tau$ of rank $n$ is a  
surjective map between two finite sets:
$$J_0\otno J_{-1}$$
equipped with a multiplicity map $\mu:J_{-1}\to \IZ_{>0}$
such that $\sum_{j\in J_{-1}}\mu(j)=n$. 
\end{defn}

The points of $J_0$ will be called {\em residual eigenspaces}, and 
the points of $J_{-1}$ will be called {\em Jordan blocks}. 
A {tame fission tree} is {\em semisimple} if 
each Jordan block has size $1$, i.e. $\mu(J_{-1})=\{1\}$.
Observe that the choice of a rank $n$ 
tame fission tree is the same as the choice of a double partition of $n$:

\begin{lem}
The choice of a  tame fission tree is the same as the choice of a 
finite set $J_0$ and a non-zero  integer partition for each $j\in J_0$.
Any matrix $A\in \GL(\IC^n)$ determines a rank $n$ tame fission tree
by taking its generalised eigenspaces, and the sizes of the Jordan blocks in each eigenspace.
More generally any automorphism of an 
$\IR$-graded vector space of dimension $n$ 
determines a rank $n$ tame fission tree
(where  the ``eigenvalues'' are now in $\IR\times \IC^*$). 
\end{lem}

\begin{defn}
An \usm{extended fission tree} $\wh \cT$ 
is the data of a fission tree $\cT$ (as in \cite{twmcg} Defn. 3.18)
together with a tame fission tree $\tau_i$ of rank $n_i$ for each 
leaf  $i\in \cT$ 
(where $n_i=\underline{n}(i)$ is the multiplicity of the leaf $i$).
An extended fission tree $\wh \cT$ is \usm{untwisted} if $\cT$ is untwisted.
It is semisimple if each $\tau_i$ is semisimple.
The slope of $\wh \cT$ is defined to be the slope 
(i.e. Poincar\'e--Katz rank $\ge 0$) of $\cT$. 
\end{defn}

\begin{figure}[ht]
    \centering
\includegraphics[width=0.6\textwidth]{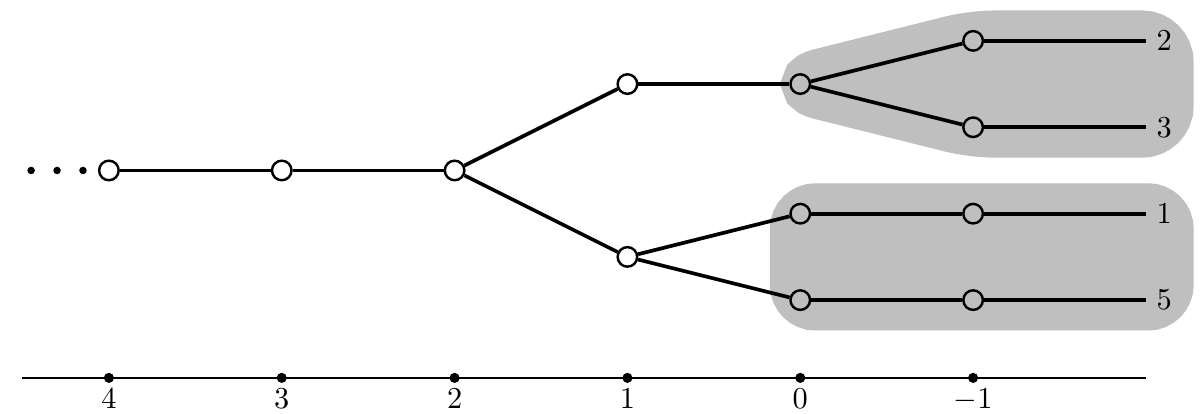}
  \caption{Extended fission tree $\widehat\IT$ 
  (slope $1$, rank $11$, \& $\IT$ has $2$ leaves).} 
    \label{fig: extfissiontree}
\end{figure}

See Fig. \ref{fig: extfissiontree} for an example. 
The untwisted extended fission trees are now easy to count, 
due to the statement:
\begin{lem}
The operation of shifting the height by $2$ gives a bijection 
between rank $n$ untwisted fission trees, and 
untwisted extended fission trees of rank $n$.
Under this correspondence the multiplicity one fission trees
correspond to the semisimple extended fission trees.
If $k\ge 1$, the extended fission trees of slope $k$ 
correspond to 
fission trees of slope $k+2$, whereas those of slope $0$ correspond 
to all the fission trees of slope $\le 2$.
\end{lem}

\begin{cor}
For any $n, k\ge 1$, the number of rank $n$ untwisted:

a)  extended fission 
trees of slope $k$  is
$\psi(k+2,n)$, and

b)  semisimple extended fission 
trees of slope $k$  is $\phi(k+2,n)$, and

c)  extended fission trees of slope $0$ is equal to
$\Phi(3,n)=\psi(2,n)+\psi(1,n)+\psi(0,n)$, and 
there are $\phi(2,n)+\phi(1,n)+\phi(0,n)=p(n)$ of them that are semisimple.
\end{cor}

Any (very) good meromorphic connection 
(in the sense of \cite{hit70, wnabh, ihptalk}) 
determines an extended fission tree at each marked point, and we expect
(cf. \cite{ballandras-mmap})  
the resulting forest of extended fission trees determines the
 deformation class of the corresponding wild nonabelian Hodge space 
 $\fM(\bSi,\wh\cC)$ .

\section{Conclusion}
In this paper we explained how to count the 
multiplicity $1$ untwisted fission trees 
by iterating the Euler transform
 and looking at the difference
of the results.
Thus the fission trees (amongst the most basic objects in 2d gauge theory)
are counted by iterating a very basic counting method  
(in effect generalising the counting of nontrivial integer partitions and double partitions to ``deeper'' cases).
This enabled the counting of several closely related objects 
(fission trees with any multiplicities, fission graphs, 
equipped fission graphs, supernova graphs, extended fission trees) and will be used in the classification of wild nonabelian Hodge moduli spaces (cf. the Lax project \cite{hit70}). 
In the sequel we will consider  automorphisms of the untwisted fission trees, in order to obtain the twisted fission trees.
One can also consider the analogous question for other 
structure groups, beyond the $\GL_n(\IC)$ case considered 
here\footnote{In effect the untwisted fission trees parameterise the possible chains of fission subgroups as in \cite{gbs} (33) p.35, and this definition is for any complex reductive group $G$ (see also \cite{wnabh} (2.2), the general discussion of fission in \cite{fission} and 
the study of 
$G$-fission trees in \cite{DRT22, doucotrembado2023topology}).}.

\appendix

\section{Squares of squares}
We saw that $\phi(3,4)=\psi(2,4)=9$ and will prove here 
that $\phi(k,4)=\psi(k-1,4)=k^2$ in general, i.e.:
\begin{prop}\label{prop: sqsq}
There are exactly $k^2$ untwisted fission trees
of rank $4$ \& slope $k-1$
(i.e. there are $k^2$  untwisted multiplicity one 
fission trees with $4$ leaves \& slope $k$).
\end{prop}
Said differently this means that:
\begin{prop}
If $k\ge 2$ then there are exactly $k^2$ fission graphs
with $4$ nodes such that the maximal 
edge multiplicity is exactly $k-1$.
\end{prop}

This will be proved by  induction. 
(In the rest of this section a fission tree means a ``multiplicity one untwisted fission tree''.)
First note two simpler facts:

\begin{lem} If $k\ge 1$ then:

$\bullet$ There is exactly $1$ 
fission tree with two leaves and slope $k$, and

$\bullet$ 
There are exactly $k$ fission trees with 
three leaves and slope $k$.
\end{lem}
\pf The first statement is clear. 
For the second: there are either one or 
two branch nodes in the 
corresponding fission tree. There is just
one with one (the generic tree).
If there are two branch nodes, 
then the only choice is 
the height of the lower branch node: 
it can be $2,3,\ldots,k$, 
so we get exactly $1+(k-1)=k$ trees.\epf

Now we can prove Proposition \ref{prop: sqsq}:

\pf
The induction we will do reduces to the identity:
$$ k^2 = (k-1)^2 + 1 + 2(k-2) +2.$$
To see this consider the map $J_2\otno J_1$ from the bottom of the fission tree.
If this map is bijective then the tree is just 
the extension 
of one of slope $k-1$ 
(increasing the length of each branch by $1$). By induction there are $(k-1)^2$ of these fission trees.

If the map $J_2\otno J_1$ is not bijective 
then $J_2$ has either  $2$ or $3$ elements.

If $J_2$ has $3$ elements then 
the tree comes from  a fission tree
of slope $k-1$ and $3$ leaves, by adding exactly one 
branch node (and extending the other branches by $1$). 
There are $1+(k-2)$ of these trees. The first 
is generic, 
and up to isomorphism 
there is just one way to extend 
it to a tree with four nodes. 
The other $k-2$ can be extended in two ways, 
so yield $2(k-2)$ trees. 

If $J_2$ has $2$ elements  
then the tree comes from the
unique fission tree
of slope $k-1$ and $2$ leaves.
There are two ways to extend it 
to get $4$ leaves ($2+2$ or $1+3$).
This gives $2$ fission trees.

This gives exactly $k^2=(k-1)^2 + 1 + 2(k-2) +2$ in total.
\epf

\section{Euler transform.}\label{sn: appx Euler transform}

For completeness we will give a proof of the Euler transform. 
This is well-known (cf. Kaneiwa \cite{kaneiwa}), 
but it might be helpful to give a
direct proof of the main facts 
in the (simple) language of multisets.

Let $S$ be a set. Recall that a ``multiset with underlying set contained in $S$'' is the same thing as a map
$$\mu:S\to \IZ_{\ge 0}$$
(the multiplicity map)
assigning an integer multiplicity $\ge 0$ to each element of $S$.
A multiset is \usm{finite} if the underlying set 
$$\supp(\mu):=\{s\in S\st \mu(s)>0\}\subset S$$
is a finite set.
The \usm{size} of a finite multiset with multiplicity map   
$\mu$ is 
$$\abs{\mu} = \sum_{s\in S}\mu(s).$$
Two multisets are equal if they have the same underlying set and the same multiplicity map (i.e. we ignore everything in $S$ with multiplicity zero).

\begin{lem}
Choose $k,n\in \IN$.
The number $\MS(k,n)$ 
of distinct multisets of size $n$ 
with
underlying set contained in $\{1,\ldots,k\}$ is given by 
$$\MS(k,n)=\bmx n+k-1 \\ n \emx=\bmx n+k-1 \\ k-1 \emx.$$
\end{lem}
\pf
Consider ordered sequences of symbols of the form:
$$
\circ \circ \circ \bigl\vert 
\circ \bigl\vert 
\bigl\vert 
\circ \circ 
\bigl\vert $$
where there are $n$ circles and $k-1$ vertical bars 
(in this example  $n=6,k=5$).
Such a sequence determines a multiset with
underlying set contained in $\{1,\ldots,k\}$. For example here we get $\mu=3,1,0,2,0$ for $1,2,3,4,5$ respectively (counting the circles from left to right).
This gives a bijection, so  $\MS(k,n)$ counts all such sequences.
But there are exactly $n+k-1$ symbols in such a sequence
and the sequence is uniquely determined by choosing which $n$
of the symbols is a circle 
(or equivalently which $k-1$ of them should be a bar). 
Thus there are $\bsmx n+k-1 \\ n \esmx=\bsmx n+k-1 \\ k-1 \esmx$
such sequences, and thus this number of multisets of the desired type.
\epf

Let $A(1),A(2),\ldots$ be a sequence of disjoint finite sets.
We will say that the elements in $A(i)$ have weight $i$, and write 
$A=\bigsqcup A(i)$. 
Thus $A$ is an $\IN$-graded set. 
Let $a(i)=\#A(i)$ be the number of elements of weight $i$.

Note that if $M$ is a finite multiset with underlying set contained in $A$ and multiplicity map $\mu:A\to \IZ_{\ge 0}$,  then 
$M$ has a well-defined weight, given by 
$$w(M)=\sum_{x\in A} \mu(x)w(x)$$
where $w(x)$ is the weight of $x\in A$.

Let $B(i)=\{M\st w(M)=i\}$ be the set of finite multisets with underlying set contained in $A$, that have weight $i$.
Define $B=\bigsqcup B(i) =:\Euler(A)$ to be the Euler transform of the graded set $A$. 
Let $b(i)=\#B(i)$ be the number of elements of $B(i)$.

\begin{lem}
$$b(i) = \sum_{\begin{smallmatrix}s_1,s_2,\ldots \ge 0\\
1\cdot s_1 + 2\cdot s_2+\cdots =i\end{smallmatrix}} 
\Prod_{k=1}^i \MS(a(k), s_k).$$
\end{lem}
\pf
This should be self-explanatory:
Any such $M$ leads to a sum of the form
$$1\cdot s_1 + 2\cdot s_2+\cdots =i$$ 
where $s_k=\#M(k)$ 
is the number of elements of $M$ of weight $k$ (counted with multiplicity).
Thus for each possible such sum we need to count the number of possible $M$'s that lead to that sum.
By definition of $\MS(k,n)$, 
the number of them is $\Prod_{k=1}^i \MS(a(k), s_k)$.
\epf

In turn the generating functions are related by the Euler transform:
\begin{prop}
$$1+\sum_1^\infty b(n)u^n = 
\Prod_1^\infty(1-u^m)^{-a(m)}\in \IZ\flb u \frb.$$
\end{prop}
\pf
This is now straightforward, using the previous two lemmas, since the binomial theorem in this context says that: 
$$(1-u^{m})^{-a} = \sum_{n=0}^\infty \MS(a,n)u^{mn}$$
for integers $a,m$. 
\epf

Thus if  $A=\bigsqcup_{i\ge 1}A(i)$ is any graded set such that 
$A(i)$ is finite, then $\Euler(A)$ is defined to be the 
set of finite multisets with underlying set contained in $A$, and there is a natural inclusion $A\into \Euler(A)$
taking $a\in A$ to the multiset given by $a$ with 
multiplicity $1$.

In effect Table 1 arises by iterating this Euler transform, starting with the set $A=\{1\}$ consisting of one element of weight $1$, and counting the new elements that arise. 
Thus 
$\Euler(A) \cong \IZ_{> 0}$
and  row 1 
of Table 1 lists the counts of elements of the set 
$\Euler(A)\setminus A$, of each weight.
In turn $\Euler^2(A)=\Euler\circ\Euler(A)$
is the set $\cP^*$ of all partitions of all integers $\ge 1$.
Row 2 
of Table  1 lists the counts of elements of the set 
$\Euler^2(A)\setminus \Euler(A)$
of each weight, i.e. the partitions with $>1$ part
($\sim$ the untwisted fission trees of slope $2$, mult. $1$, weighted by their leaf counts).
Similarly row k 
of Table  1 lists the counts of elements of the set 
$\Euler^k(A)\setminus \Euler^{k-1}(A)$, i.e. 
the (multiplicity $1$) untwisted fission trees of slope $k$.

\renewcommand{\baselinestretch}{1}              %
\normalsize
\bibliographystyle{amsplain}    \label{biby}
\bibliography{../thesis/syr}

\vspace{0.1cm} 

\noindent

\ 

\noindent
{\small
\noindent
Universit\'e Paris Cit\'e and Sorbonne Universit\'e,
CNRS, IMJ-PRG, F-75013 Paris, France} \\
\noindent
boalch@imj-prg.fr \hfill
\url{https://webusers.imj-prg.fr/~philip.boalch/}

\end{document}